\documentclass[12pt,reqno]{article}

\usepackage{amsmath,amssymb,amsthm}
\usepackage[all]{xy}

\usepackage{hyperref}

\usepackage{fullpage}
\usepackage{float}

\newcommand{\seqnum}[1]{\href{http://www.research.att.com/cgi-bin/access.cgi/as/~njas/sequences/eisA.cgi?Anum=#1}{\underline{#1}}}

\newtheorem{thm}{Theorem}[section]
\newtheorem{prop}[thm]{Proposition}
\newtheorem{lem}[thm]{Lemma}
\newtheorem{cor}[thm]{Corollary}
\newtheorem{conj}[thm]{Conjecture}
\theoremstyle{definition}
\newtheorem{defn}[thm]{Definition}
\newtheorem{ex}[thm]{Example}
\newtheorem{rem}[thm]{Remark}
\newcommand{\C}{\mathcal{C}}
\renewcommand{\H}{\mathfrak{H}}
\newcommand{\Q}{\mathbb{Q}}
\newcommand{\R}{\mathbb{R}}
\newcommand{\Z}{\mathbb{Z}}
\newcommand{\ZZ}{\mathcal{Z}}
\DeclareMathOperator{\spa}{span}
\DeclareMathOperator{\Ker}{Ker}
\DeclareMathOperator{\sgn}{sgn}
\renewcommand{\Im}{\operatorname{Im}}
\DeclareMathOperator{\rank}{rank}
\newcommand{\ve}[1]{\boldsymbol{#1}}

\begin{document}

\begin{center}
\vskip 1cm{\LARGE\textbf{Combinatorial remarks on the cyclic sum formula for multiple zeta values}}
\vskip 1cm
\large

Shingo Saito\\
Faculty of Mathematics\\
Kyushu University\\
744, Motooka, Nishi-ku, Fukuoka, 819--0395, Japan\\
\href{mailto:ssaito@math.kyushu-u.ac.jp}{\texttt{ssaito@math.kyushu-u.ac.jp}}\\
\ \\
Tatsushi Tanaka\\
Faculty of Mathematics\\
Kyushu University\\
744, Motooka, Nishi-ku, Fukuoka, 819--0395, Japan\\
\href{mailto:t.tanaka@math.kyushu-u.ac.jp}{\texttt{t.tanaka@math.kyushu-u.ac.jp}}\\
\ \\
Noriko Wakabayashi\\
Faculty of Engineering\\
Kyushu Sangyo University\\
3--1, Matsukadai 2-chome, Higashi-ku, Fukuoka, 813--8503, Japan\\
\href{mailto:noriko@ip.kyusan-u.ac.jp}{\texttt{noriko@ip.kyusan-u.ac.jp}}\\
\end{center}

\vskip .2 in

\begin{abstract}
 The multiple zeta values are generalizations of the values of the Riemann zeta function
 at positive integers.
 They are known to satisfy a number of relations,
 among which are the cyclic sum formula.
 The cyclic sum formula can be stratified via linear operators
 defined by the second and third authors.
 We give the number of relations belonging to each stratum
 by combinatorial arguments.
\end{abstract}

\renewcommand{\theequation}{\arabic{equation}}

\section{Introduction}
The Riemann zeta function
\[
 \zeta(s)=\sum_{n=1}^{\infty}\frac{1}{n^s}
\]
is one of the most important functions in mathematics,
and its values $\zeta(k)$ at $k\in\Z_{\ge2}$ are among the most interesting real numbers.
The \emph{multiple zeta values} (MZVs for short) are defined as generalizations of the values $\zeta(k)$:
\begin{defn}
 For $k_1,\ldots,k_l\in\Z_{\ge1}$ with $k_1\ge2$,
 the \emph{multiple zeta value} $\zeta(k_1,\ldots,k_l)$ is a real number defined by
 \[
  \zeta(k_1,\ldots,k_l)=\sum_{n_1>\cdots >n_l\ge1}\frac{1}{n_1^{k_1}\cdots n_l^{k_l}}.
 \]
\end{defn}

Little is known about the irrationality or transcendentality of MZVs;
we however know that MZVs abound with relations among themselves,
the simplest example being $\zeta(2,1)=\zeta(3)$.

The \emph{cyclic sum formula} (CSF for short), described in Section~\ref{sec:CSF},
is a class of $\Q$-linear relations among MZVs, established by Hoffman-Ohno~\cite{HO}.
Their proof appealed to partial fraction decomposition,
whereas another proof given by the second and third authors~\cite{TW} of the present article
proceeded by showing that the CSF is included in Kawashima's relation,
which is believed to be rich enough to yield all relations among MZVs.

This present paper is aimed at providing combinatorial arguments to find the ranks
of linear operators defined in \cite{TW}.
In order to facilitate access for both algebraists and combinatorists,
we strive to make the exposition as self-contained as possible.

\subsection*{Sets of multi-indices}
The study of MZVs inevitably requires frequent use of multi-indices.
We here summarize the sets of multi-indices used in this paper.
For $k,l\in\Z_{\ge1}$, put
\begin{align*}
 I_{k,l}^1&=\{(k_1,\ldots,k_l)\in\Z_{\ge1}^l\mid k_1+\cdots+k_l=k\},\\
 I_{k,l}^0&=\{(k_1,\ldots,k_l)\in I_{k,l}^1\mid k_1\ge2\},\\
 \check{I}_{k,l}^1&=\{(k_1,\ldots,k_l)\in I_{k,l}^1\mid\text{not all of $k_1,\ldots,k_l$ are $1$}\};
\end{align*}
for $k\in\Z_{\ge1}$, put
\[
 I_k^1=\bigcup_{l=1}^{\infty}I_{k,l}^1,\qquad
 I_k^0=\bigcup_{l=1}^{\infty}I_{k,l}^0,\qquad
 \check{I}_k^1=\bigcup_{l=1}^{\infty}\check{I}_{k,l}^1;
\]
put
\[
 I^1=\bigcup_{k,l=1}^{\infty}I_{k,l}^1,\qquad
 I^0=\bigcup_{k,l=1}^{\infty}I_{k,l}^0,\qquad
 \check{I}^1=\bigcup_{k,l=1}^{\infty}\check{I}_{k,l}^1.
\]

The elements of $I_{k,l}^1$ are said to have \emph{weight} $k$ and \emph{depth} $l$.
For each $\ve{k}=(k_1,\ldots,k_l)\in I^0$
the MZV $\zeta(\ve{k})=\zeta(k_1,\ldots,k_l)$ is defined.

\section{Hoffman's algebra}
In the discussion of MZVs, it is convenient to use the algebra introduced by Hoffman.
Let $\H=\Q\langle x,y \rangle$ be the noncommutative polynomial algebra over $\Q$
in two indeterminates $x$ and $y$,
and put $\H^0=\Q+x\H y$, which is a subalgebra of $\H$.
Write $z_k=x^{k-1}y\in\H$ for $k\in\Z_{\ge1}$, so that
\[
 \{1\}\cup\{z_{k_1}\cdots z_{k_l}\mid \ve{k}=(k_1,\ldots,k_l)\in I^0\}
\]
is a $\Q$-vector space basis for $\H^0$.
It follows that we may define a $\Q$-linear map $Z\colon\H^0\to\R$ by
setting $Z(1)=1$ and $Z(z_{k_1}\cdots z_{k_l})=\zeta(\ve{k})$
for $\ve{k}=(k_1,\ldots,k_l)\in I^0$.

The MZVs are known to fulfill a number of $\Q$-linear relations,
each of which corresponds to an element of $\Ker Z\subset\H^0$.
Since Goncharov~\cite{Go} conjectured that
the MZVs $\zeta(\ve{k})$ with $\ve{k}$ having different weights
are $\Q$-linearly independent,
we look at the $\Q$-vector subspaces defined by
\begin{align*}
 \H_k^0&=\spa_{\Q}\{z_{k_1}\cdots z_{k_l}\mid \ve{k}=(k_1,\ldots,k_l)\in I_k^0\}\\
 &=\{w\in\H^0\mid\text{$w$ is a homogeneous polynomial of degree $k$}\}\cup\{0\}\subset\H^0,\\
 \ZZ_k&=Z(\H_k^0)=\spa_{\Q}\{\zeta(\ve{k})\mid \ve{k}\in I_k^0\}\subset\R
\end{align*}
for each $k\in\Z_{\ge1}$.

Let $(d_k)_{k\ge1}$ be the Padovan sequence (\seqnum{A000931})
defined by $d_1=0$, $d_2=d_3=1$, and $d_k=d_{k-2}+d_{k-3}$ for $k\ge4$.
Zagier~\cite{Z} conjectured the following:
\begin{conj}[Zagier~\cite{Z}]\label{conj:Z}
 We have $\dim\ZZ_k=d_k$ for all $k\in\Z_{\ge1}$.
\end{conj}

Let $k\in\Z_{\ge2}$.
In light of the fact that $\dim\H_k^0=\#I_k^0=2^{k-2}$,
Conjecture~\ref{conj:Z} means that the MZVs must satisfy plenty of $\Q$-linear relations.
Note that Conjecture~\ref{conj:Z} is equivalent to saying that
the restriction $Z\vert_{\H_k^0}\colon\H_k^0\to\R$ of $Z$ to $\H_k^0$ has rank $d_k$,
which is also equivalent to
\[
 \dim\Ker Z\vert_{\H_k^0}=\dim(\Ker Z\cap\H_k^0)=2^{k-2}-d_k.
\]

\begin{table}[h]
 \caption{Dimensions concerning $Z\vert_{\H_k^0}$}
 \begin{center}
  \begin{tabular}{c|ccccccccc|c}
   $k$&$2$&$3$&$4$&$5$&$6$&$7$&$8$&$9$&$10$&Sequence Number\\\hline
   $2^{k-2}(=\dim\H_k^0)$&$1$&$2$&$4$&$8$&$16$&$32$&$64$&$128$&$256$&\seqnum{A000079}\\
   $d_k(\stackrel{?}{=}\rank Z\vert_{\H_k^0})$&$1$&$1$&$1$&$2$&$2$&$3$&$4$&$5$&$7$&\seqnum{A000931}\\
   $2^{k-2}-d_k(\stackrel{?}{=}\dim\Ker Z\vert_{\H_k^0})$&$0$&$1$&$3$&$6$&$14$&$29$&$60$&$123$&$249$&\seqnum{A038360}
  \end{tabular}
 \end{center}
\end{table}

Goncharov~\cite{G} and Terasoma~\cite{T} partially proved Conjecture~\ref{conj:Z}:
\begin{thm}[Goncharov~\cite{G}, Terasoma~\cite{T}]\label{thm:GT}
 We have $\dim\ZZ_k\le d_k$ for all $k\in\Z_{\ge1}$.
\end{thm}

Since their proofs of Theorem~\ref{thm:GT} resort to algebraic geometry
and fails to give concrete relations among MZVs,
it still lies at the heart of research to find sufficiently many $\Q$-linear relations.
Also, the converse inequality is far from being solved.

\section{Cyclic sum formula}\label{sec:CSF}
Numerous concrete $\Q$-linear relations among MZVs have been obtained so far,
and our focus is on the following \emph{cyclic sum formula} (CSF for short),
first proved by Hoffman-Ohno~\cite{HO}:
\begin{thm}[Cyclic sum formula]\label{thm:CSF}
 If $(k_1,\ldots,k_l)\in\check{I}^1$, then
 \[
  \sum_{j=1}^{l}\sum_{i=1}^{k_j-1}\zeta(k_j-i+1,k_{j+1},\ldots,k_l,k_1,\dots,k_{j-1},i)
  =\sum_{j=1}^{l}\zeta(k_j+1,k_{j+1},\ldots,k_l,k_1,\ldots,k_{j-1}).
 \]
\end{thm}

\begin{ex}
 The cyclic sum formula for $l=1$ and $k_1=2$ gives
 $\zeta(2,1)=\zeta(3)$.
\end{ex}

In dealing with the CSF, it is convenient to extend the indices of $k_j$ to all $j\in\Z$
by declaring $k_j=k_{j'}$ whenever $j\equiv j'\pmod{l}$.
Then the CSF can simply be written as
\[
 \sum_{j=1}^{l}\sum_{i=1}^{k_j-1}\zeta(k_j-i+1,k_{j+1},\ldots,k_{j+l-1},i)
 =\sum_{j=1}^{l}\zeta(k_j+1,k_{j+1},\ldots,k_{j+l-1}).
\]
This convention will be used tacitly throughout the paper.

In order to describe the CSF in terms of Hoffman's algebra, we write
\begin{align*}
 \check{\H}^1&=\spa_{\Q}\{z_{k_1}\cdots z_{k_l}\mid \ve{k}=(k_1,\ldots,k_l)\in\check{I}^1\}\\
 &=\spa_{\Q}\{w\in\H\mid\text{$w$ is a monomial ending with $y$ but not a power of $y$}\}\subset\H
\end{align*}
and define a $\Q$-linear map $\rho\colon\check{\H}^1\to x\H y\subset\H^0$ by setting
\[
 \rho(z_{k_1}\cdots z_{k_l})
 =\sum_{j=1}^{l}\sum_{i=1}^{k_j-1}z_{k_j-i+1}z_{k_{j+1}}\cdots z_{k_{j+l-1}}z_i
 -\sum_{j=1}^{l}z_{k_j+1}z_{k_{j+1}}\cdots z_{k_{j+l-1}}
\]
for $\ve{k}=(k_1,\ldots,k_l)\in\check{I}^1$.
Then the CSF is equivalent to saying that
$\Im\rho\subset\Ker Z$.

For each $k\in\Z_{\ge1}$, if we put
\begin{align*}
 \check{\H}_k^1&=\spa_{\Q}\{z_{k_1}\cdots z_{k_l}\mid(k_1,\ldots,k_l)\in\check{I}_k^1\}\\
 &=\{w\in\check{\H}^1\mid\text{$w$ is a homogeneous polynomial of degree $k$}\}\cup\{0\},
\end{align*}
then $\rho$ satisfies that $\rho(\check{\H}_k^1)\subset\H_{k+1}^0$.
Therefore, it follows from Theorem~\ref{thm:CSF} that
\[
 \rho(\check{\H}_{k-1}^1)\subset\Ker Z\cap\H_k^0
\]
for all $k\in\Z_{\ge2}$.

In view of Conjecture~\ref{conj:Z},
which is equivalent to $\dim(\Ker Z\cap\H_k^0)=2^{k-2}-d_k$ for $k\in\Z_{\ge2}$,
it is natural to ask for $\dim\rho(\check{\H}_{k-1}^1)$
because it can be regarded as the number of relations given by the CSF.
The following theorem is a folklore:
\begin{thm}\label{thm:dim_CSF}
 For $k\in\Z_{\ge2}$, we have
 \[
  \dim\rho(\check{\H}_{k-1}^1)=\frac{1}{k-1}\sum_{m\mid k-1}\varphi\biggl(\frac{k-1}{m}\biggr)2^m-2,
 \]
 where $\varphi$ denotes Euler's totient function.
\end{thm}

\begin{table}[h]
 \caption{Dimensions concerning the CSF}
 \begin{center}
  \begin{tabular}{c|ccccccccc|c}
   $k$&$2$&$3$&$4$&$5$&$6$&$7$&$8$&$9$&$10$&Sequence Number\\\hline
   $2^{k-2}-d_k\bigl(\stackrel{?}{=}\dim(\Ker Z\cap\H_k^0)\bigr)$&%
   $0$&$1$&$3$&$6$&$14$&$29$&$60$&$123$&$249$&\seqnum{A038360}\\
   $\dim\rho(\check{\H}_{k-1}^1)$&$0$&$1$&$2$&$4$&$6$&$12$&$18$&$34$&$58$&\seqnum{A052823}
  \end{tabular}
 \end{center}
\end{table}

\section{The operator $\rho_n$ and the statement of our main theorem}
\subsection{The operator $\rho_n$}
The second and third authors~\cite{TW} of the present article
defined linear maps $\rho_n\colon\H\to\H$
with the aim of giving an algebraic proof of the CSF
by reducing it to Kawashima's relation.
We will not elaborate on their proof here,
but focus on the stratification of the CSF provided by $\rho_n$.
Note that our usage of indices is different from that in~\cite{TW}:
what we mean by $\rho_n$ is denoted by $\rho_{n+1}$ in \cite{TW}.

Let $n\in\Z_{\ge0}$ and consider the $(n+2)$nd tensor power $\H^{\otimes(n+2)}$
of $\H$ over $\Q$.
We first make $\H^{\otimes(n+2)}$ an $\H$-bimodule by setting
\[
 a\diamond(w_1\otimes\dots\otimes w_{n+2})\diamond b
 =w_1b\otimes w_2\otimes\dots\otimes w_{n+1}\otimes aw_{n+2}
\]
for $a,b,w_1,\dots,w_{n+2}\in\H$.
Writing $z=x+y$, we define a $\Q$-linear map $\C_n\colon\H\to\H^{\otimes(n+2)}$ by
setting $\C_n(1)=0$, $\C_n(x)=x\otimes z^{\otimes n}\otimes y$,
$\C_n(y)=-x\otimes z^{\otimes n}\otimes y$, and
\[
 \C_n(ww')=\C_n(w)\diamond w'+w\diamond\C_n(w').
\]
We next define a $\Q$-linear map $M_n\colon\H^{\otimes(n+2)}\to\H$ by setting
\[
 M_n(w_1\otimes\cdots\otimes w_{n+2})=w_1\cdots w_{n+2}.
\]
Finally, set $\rho_n=M_n\circ\C_n\colon\H\to\H$.

\begin{rem}
 The recurrence relation in the definition of $\C_n$ shows that
 \[
  \C_n(w_1\cdots w_k)=
  \sum_{j=1}^{k}\bigl(w_1\cdots w_{j-1}\diamond\C_n(w_j)\diamond w_{j+1}\cdots w_k\bigr)
 \]
 for $w_1,\ldots,w_k\in\H$.
\end{rem}

\begin{prop}\label{prop:rho0=rho}
 We have $\rho_0(w)=\rho(w)$ for all $w\in\check{\H}^1$.
\end{prop}

\begin{proof}
 We may assume that $w=z_{k_1}\cdots z_{k_l}$ for some
 $\ve{k}=(k_1,\ldots,k_l)\in\check{I}^1$.
 For $k\in\Z_{\ge1}$, we have
 \begin{align*}
  \C_0(z_k)&=\C_0(x^{k-1}y)
  =\sum_{i=1}^{k-1}\bigl(x^{i-1}\diamond\C_0(x)\diamond x^{k-i-1}y\bigr)+x^{k-1}\diamond\C_0(y)\\
  &=\sum_{i=1}^{k-1}\bigl(x^{i-1}\diamond(x\otimes y)\diamond x^{k-i-1}y\bigr)
    +x^{k-1}\diamond(-x\otimes y)\\
  &=\sum_{i=1}^{k-1}(x^{k-i}y\otimes x^{i-1}y)-x\otimes x^{k-1}y
   =\sum_{i=1}^{k-1}(z_{k-i+1}\otimes z_i)-x\otimes z_k.
 \end{align*}
 It follows that
 \begin{align*}
  \C_0(w)
  &=\C_0(z_{k_1}\cdots z_{k_l})
   =\sum_{j=1}^{l}\bigl(z_{k_1}\cdots z_{k_{j-1}}\diamond\C_0(z_{k_j})\diamond z_{k_{j+1}}\cdots z_{k_l}\bigr)\\
  &=\sum_{j=1}^{l}\sum_{i=1}^{k_j-1}
    \bigl(z_{k_1}\cdots z_{k_{j-1}}\diamond(z_{k_j-i+1}\otimes z_i)\diamond z_{k_{j+1}}\cdots z_{k_l}\bigr)\\
  &\qquad-\sum_{j=1}^{l}
   \bigl(z_{k_1}\cdots z_{k_{j-1}}\diamond(x\otimes z_{k_j})\diamond z_{k_{j+1}}\cdots z_{k_l}\bigr)\\
  &=\sum_{j=1}^{l}\sum_{i=1}^{k_j-1}
    (z_{k_j-i+1}z_{k_{j+1}}\cdots z_{k_l}\otimes z_{k_1}\cdots z_{k_{j-1}}z_{i})
  -\sum_{j=1}^{l}(xz_{k_{j+1}}\cdots z_{k_l}\otimes z_{k_1}\cdots z_{k_j}),
 \end{align*}
 and so
 \begin{align*}
  \rho_0(w)
  &=M_0\bigl(\C_0(w)\bigr)\\
  &=\sum_{j=1}^{l}\sum_{i=1}^{k_j-1}
   z_{k_j-i+1}z_{k_{j+1}}\cdots z_{k_l}z_{k_1}\cdots z_{k_{j-1}}z_i
   -\sum_{j=1}^{l}xz_{k_{j+1}}\cdots z_{k_l}z_{k_1}\cdots z_{k_j}\\
  &=\sum_{j=1}^{l}\sum_{i=1}^{k_j-1}
   z_{k_j-i+1}z_{k_{j+1}}\cdots z_{k_{j+l-1}}z_i
   -\sum_{j=1}^{l}xz_{k_{j+1}}\cdots z_{k_{j+l}},
 \end{align*}
 which is equal to $\rho(w)$ because
 \[
  \sum_{j=1}^{l}xz_{k_{j+1}}\cdots z_{k_{j+l}}
  =\sum_{j=1}^{l}xz_{k_j}\cdots z_{k_{j+l-1}}
  =\sum_{j=1}^{l}z_{k_j+1}z_{k_{j+1}}\cdots z_{k_{j+l-1}}.\qedhere
 \]
\end{proof}

\subsection{Properties of $\rho_n$}
\begin{defn}
 Define $\sgn\colon\{x,y,z\}\to\{1,-1,0\}$ by
 \[
  \sgn(u)=
  \begin{cases}
   1&\text{if $u=x$};\\
   -1&\text{if $u=y$};\\
   0&\text{if $u=z$},
  \end{cases}
 \]
 so that $\C_n(u)=\sgn(u)(x\otimes z^{\otimes n}\otimes y)$
 for $u\in\{x,y,z\}$.
\end{defn}

\begin{lem}\label{lem:rho_n_prod_xyz}
 If $w=u_1\cdots u_k$, where $u_1,\ldots,u_k\in\{x,y,z\}$, then
 \[
  \rho_n(w)=\sum_{j=1}^{k}\sgn(u_j)xu_{j+1}\cdots u_kz^nu_1\cdots u_{j-1}y.
 \]
\end{lem}

\begin{proof}
 We have
 \begin{align*}
  \C_n(w)&=\C_n(u_1\cdots u_k)
  =\sum_{j=1}^{k}\bigl(u_1\cdots u_{j-1}\diamond\C_n(u_j)\diamond u_{j+1}\cdots u_k\bigr)\\
  &=\sum_{j=1}^{k}\sgn(u_j)
    \bigl(u_1\cdots u_{j-1}\diamond(x\otimes z^{\otimes n}\otimes y)\diamond u_{j+1}\cdots u_k\bigr)\\
  &=\sum_{j=1}^{k}\sgn(u_j)
    (xu_{j+1}\cdots u_k\otimes z^{\otimes n}\otimes u_1\cdots u_{j-1}y),
 \end{align*}
 and so
 \[
  \rho_n(w)=M_n\bigl(\C_n(w)\bigr)
  =\sum_{j=1}^{k}\sgn(u_j)xu_{j+1}\cdots u_kz^nu_1\cdots u_{j-1}y.\qedhere
 \]
\end{proof}

\begin{prop}\label{prop:rho_n_n+1}
 We have
 \[
  \rho_{n+1}(w)=\rho_n(zw)
 \]
 for all $w\in\H$.
\end{prop}

\begin{proof}
 We may assume that $w=u_1\cdots u_k$ for some $u_1,\ldots,u_k\in\{x,y\}$.
 Then Lemma~\ref{lem:rho_n_prod_xyz} shows that
 \begin{align*}
  \rho_n(zw)
  &=\sgn(z)xu_1\cdots u_kz^ny+\sum_{j=1}^{k}\sgn(u_j)xu_{j+1}\cdots u_kz^nzu_1\cdots u_{j-1}y\\
  &=\sum_{j=1}^{k}\sgn(u_j)xu_{j+1}\cdots u_kz^{n+1}u_1\cdots u_{j-1}y
   =\rho_{n+1}(w).\qedhere
 \end{align*}
\end{proof}

\begin{cor}\label{cor:stratification}
 We have
 \[
  \{0\}=\rho_{k-2}(\check{\H}_1^1)\subset\rho_{k-3}(\check{\H}_2^1)
  \subset\cdots\subset\rho_0(\check{\H}_{k-1}^1)=\rho(\check{\H}_{k-1}^1).
 \]
 for $k\in\Z_{\ge2}$.
\end{cor}

\begin{proof}
 For each $n\in\{0,\ldots,k-3\}$,
 we have $\rho_{n+1}(\check{\H}_{k-n-2}^1)\subset\rho_n(\check{\H}_{k-n-1}^1)$
 by Proposition~\ref{prop:rho_n_n+1} and by the fact that
 if $w\in\check{\H}_{k-n-2}^1$, then $zw\in\check{\H}_{k-n-1}^1$.
\end{proof}

\subsection{Statement of our main theorem}
Corollary~\ref{cor:stratification} can be interpreted as stratifying
the $\Q$-linear relations provided by the CSF.
Since Theorem~\ref{thm:dim_CSF} tells us the dimension of the whole space,
we may well wish to find the dimensions of the subspaces $\rho_n(\check{\H}_k^1)$ in general.
Our main theorem (Theorem~\ref{thm:main}) provides a complete solution to this problem,
and it uses the following generalization of the Lucas sequence:

\begin{defn}
 For $n\in\Z_{\ge1}$, the \emph{$n$-step Lucas sequence} $(L_m^{n})_{m\ge1}$
 is defined by
 \[
  L_m^{n}=
  \begin{cases}
   2^m-1,&\text{for $m=1,\ldots,n$};\\
   L_{m-1}^{n}+\cdots+L_{m-n}^{n},&\text{for $m\ge n+1$}.
  \end{cases}
 \]
 We adopt the convention that $L_m^0=0$ for all $m\in\Z_{\ge1}$.
\end{defn}

\begin{table}[h]
 \caption{$n$-step Lucas sequences}
 \begin{center}
  \begin{tabular}{c|ccccccc|c}
   $n$&$1$&$2$&$3$&$4$&$5$&$6$&$7$&Sequence Number\\\hline
   $L_m^{1}$&$1$&$1$&$1$&$1$&$1$&$1$&$1$&\seqnum{A000012}\\
   $L_m^{2}$&$1$&$3$&$4$&$7$&$11$&$18$&$29$&\seqnum{A000032}, \seqnum{A000204}\\
   $L_m^{3}$&$1$&$3$&$7$&$11$&$21$&$39$&$71$&\seqnum{A001644}\\
   $L_m^{4}$&$1$&$3$&$7$&$15$&$26$&$51$&$99$&\seqnum{A073817}, \seqnum{A001648}\\
   $L_m^{5}$&$1$&$3$&$7$&$15$&$31$&$57$&$113$&\seqnum{A074048}, \seqnum{A023424}
  \end{tabular}
 \end{center}
\end{table}

\begin{lem}\label{lem:n-step}
 Let $n\in\Z_{\ge0}$.
 \begin{enumerate}
  \item We have $L_{n+1}^{n}=2^{n+1}-n-2$.
  \item We have $L_m^{n}=2L_{m-1}^{n}-L_{m-n-1}^{n}$ for $m\ge n+2$.
 \end{enumerate}
\end{lem}

\begin{proof}
 Easy.
\end{proof}

\begin{thm}[Main Theorem]\label{thm:main}
 If $n\in\Z_{\ge0}$ and $k\in\Z_{\ge1}$, then
 \[
  \dim\rho_n(\check{\H}_{k}^1)
  =\frac{1}{n+k}\sum_{m\mid n+k}\varphi\biggl(\frac{n+k}{m}\biggr)(2^m-L_m^n)-2.
 \]
\end{thm}

\begin{rem}
 We may easily see that this theorem is a generalization of Theorem~\ref{thm:dim_CSF}.
\end{rem}

\section{Proof of our main theorem}
\subsection{Cyclic equivalence and $\rho_0$}
\begin{defn}
 Let $l\in\Z_{\ge1}$.
 Set $X_l=\{y,z\}^l$, and let the cyclic group $\Z/l\Z$ of order $l$
 act on $X_l$ by cyclic shifts, i.e.,
 \[
  j(u_1,\ldots,u_l)=(u_{j+1},\ldots,u_{j+l})
 \]
 for $j\in\Z/l\Z$ and $(u_1,\ldots,u_l)\in X_l$.
 The equivalence relation on $X_l$ induced by the action is called the \emph{cyclic equivalence}
 and denoted by $\sim$.
 Put $Y_l=X_l/{\sim}$.
\end{defn}

\begin{prop}
 For each $l\in\Z_{\ge1}$,
 we may define $\tilde{\rho}_0\colon Y_l\to\H$ by setting
 \[
  \tilde{\rho}_0\bigl([(u_1,\ldots,u_l)]\bigr)=\rho_0(u_1\cdots u_l)
 \]
 for $(u_1,\ldots,u_l)\in X_l$.
\end{prop}

\begin{proof}
 Lemma~\ref{lem:rho_n_prod_xyz} shows that if $(u_1,\ldots,u_l)\in X_l$, then
 \begin{align*}
  \rho_0(u_1\cdots u_l)
  &=\sum_{j=1}^{l}\sgn(u_j)xu_{j+1}\cdots u_lu_1\cdots u_{j-1}y\\
  &=\sum_{\substack{j=1,\ldots,l\\(v_1,\ldots,v_l)=j(u_1,\ldots,u_l)}}\sgn(v_1)xv_2\cdots v_ly.
 \end{align*}
 It follows that $\rho_0(u_1\cdots u_l)$ depends only on
 the equivalence class $[(u_1,\ldots,u_l)]\in Y_l$,
 as required.
\end{proof}

\begin{ex}
 For $l=4$, we have
 \begin{align*}
  \tilde{\rho}_0\bigl([(z,z,z,z)]\bigr)&=0,\\
  \tilde{\rho}_0\bigl([(z,z,z,y)]\bigr)&=-xzzzy,\\
  \tilde{\rho}_0\bigl([(z,z,y,y)]\bigr)&=-x(yzz+zzy)y,\\
  \tilde{\rho}_0\bigl([(z,y,z,y)]\bigr)&=-x(zyz+zyz)y,\\
  \tilde{\rho}_0\bigl([(z,y,y,y)]\bigr)&=-x(yyz+yzy+zyy)y,\\
  \tilde{\rho}_0\bigl([(y,y,y,y)]\bigr)&=-x(yyy+yyy+yyy+yyy)y.
 \end{align*}
\end{ex}

\begin{prop}\label{prop:rho_U_lin_indep}
 The family
 $\bigl\{\tilde{\rho}_0(U)\bigm|U\in Y_l\setminus\{[(\underbrace{z,\ldots,z}_{l})]\}\bigr\}$
 is $\Q$-linearly independent.
\end{prop}

\begin{proof}
 This is because the values $\tilde{\rho}_0(U)$
 for different equivalence classes $U\in Y_l\setminus\{[(z,\ldots,z)]\}$
 are nonzero and consist of different monomials.
\end{proof}

\subsection{Relationship between $\rho_n$ and $\rho_0$}
\begin{defn}
 For $l\in\Z_{\ge1}$ and $n\in\{0,\ldots,l\}$,
 we write $X_{l,n}$ for the subset of $X_l$ consisting of
 all $(u_1,\ldots,u_l)\in X_l$ that contain at least $n$ consecutive $z$'s
 when written cyclically.
\end{defn}

\begin{ex}
 For $l=4$, we have
 \begin{align*}
  X_{4,0}&=\{(z,z,z,z),\\
   &\phantom{{}=\{}(z,z,z,y),(z,z,y,z),(z,y,z,z),(y,z,z,z),\\
   &\phantom{{}=\{}(z,z,y,y),(z,y,y,z),(y,y,z,z),(y,z,z,y),\\
   &\phantom{{}=\{}(z,y,z,y),(y,z,y,z),\\
   &\phantom{{}=\{}(z,y,y,y),(y,y,y,z),(y,y,z,y),(y,z,y,y),\\
   &\phantom{{}=\{}(y,y,y,y)\}\\
   &=[(z,z,z,z)]\cup[(z,z,z,y)]\cup[(z,z,y,y)]\cup[(z,y,z,y)]\cup[(z,y,y,y)]\cup[(y,y,y,y)]\\
   &=X_4,\\
  X_{4,1}&=\{(z,z,z,z),\\
   &\phantom{{}=\{}(z,z,z,y),(z,z,y,z),(z,y,z,z),(y,z,z,z),\\
   &\phantom{{}=\{}(z,z,y,y),(z,y,y,z),(y,y,z,z),(y,z,z,y),\\
   &\phantom{{}=\{}(z,y,z,y),(y,z,y,z),\\
   &\phantom{{}=\{}(z,y,y,y),(y,y,y,z),(y,y,z,y),(y,z,y,y)\}\\
   &=[(z,z,z,z)]\cup[(z,z,z,y)]\cup[(z,z,y,y)]\cup[(z,y,z,y)]\cup[(z,y,y,y)],\\
  X_{4,2}&=\{(z,z,z,z),\\
   &\phantom{{}=\{}(z,z,z,y),(z,z,y,z),(z,y,z,z),(y,z,z,z),\\
   &\phantom{{}=\{}(z,z,y,y),(z,y,y,z),(y,y,z,z),(y,z,z,y)\}\\
   &=[(z,z,z,z)]\cup[(z,z,z,y)]\cup[(z,z,y,y)],\\
  X_{4,3}&=\{(z,z,z,z),\\
   &\phantom{{}=\{}(z,z,z,y),(z,z,y,z),(z,y,z,z),(y,z,z,z)\}\\
   &=[(z,z,z,z)]\cup[(z,z,z,y)],\\
  X_{4,4}&=\{(z,z,z,z)\}\\
   &=[(z,z,z,z)].
 \end{align*}
\end{ex}

\begin{rem}
 Each $X_{l,n}$ is invariant under the action of $\Z/l\Z$,
 which allows us to make the following definition.
\end{rem}

\begin{defn}
 For $l\in\Z_{\ge1}$ and $n\in\{0,\ldots,l\}$, we write
 \[
  Y_{l,n}=X_{l,n}/{\sim}
  =\bigl\{[(\underbrace{z,\ldots,z}_{n},u_{n+1},\ldots,u_l)]\in Y_l
   \bigm|u_{n+1},\ldots,u_l\in\{y,z\}\bigr\}.
 \]
\end{defn}

\begin{prop}\label{prop:dim_and_Y}
 For $n\in\Z_{\ge0}$ and $k\in\Z_{\ge1}$, we have
 \[
  \dim\rho_n(\check{\H}_k^1)=\#Y_{n+k,n}-2.
 \]
\end{prop}

\begin{proof}
 Since
 \begin{align*}
  \check{\H}_k^1
  &=\spa_{\Q}\{z_{k_1}\cdots z_{k_l}\mid(k_1,\ldots,k_l)\in\check{I}_k^1\}\\
  &=\spa_{\Q}\bigl\{u_1\cdots u_{k-1}y\bigm|
    (u_1,\ldots,u_{k-1})\in\{x,y\}^{k-1}\setminus\{(\underbrace{y,\ldots,y}_{k-1})\}\bigr\}\\
  &=\spa_{\Q}\bigl\{u_1\cdots u_{k-1}y-y^k\bigm|
    (u_1,\ldots,u_{k-1})\in\{y,z\}^{k-1}\bigr\},
 \end{align*}
 Proposition~\ref{prop:rho_n_n+1} shows that
 \begin{align*}
  \rho_n(\check{\H}_k^1)
  &=\spa_{\Q}\bigl\{\rho_n(u_1\cdots u_{k-1}y)-\rho_n(y^k)\bigm|
    (u_1,\ldots,u_{k-1})\in\{y,z\}^{k-1}\bigr\}\\
  &=\spa_{\Q}\bigl\{\rho_0(z^nu_1\cdots u_{k-1}y)-\rho_0(z^ny^k)\bigm|
    (u_1,\ldots,u_{k-1})\in\{y,z\}^{k-1}\bigr\}\\
  &=\spa_{\Q}\bigl\{\tilde{\rho}_0(U)-
    \tilde{\rho}_0\bigl([(\underbrace{z,\ldots,z}_{n},\underbrace{y,\ldots,y}_k)]\bigr)\bigm|
    U\in Y_{n+k,n}\setminus\{[(\underbrace{z,\ldots,z}_{n+k})]\}\bigr\},
 \end{align*}
 which implies that $\dim\rho_n(\check{\H}_k^1)=\#Y_{n+k,n}-2$
 because of Proposition~\ref{prop:rho_U_lin_indep}.
\end{proof}

\subsection{Calculation of $\#Y_{l,n}$}
Proposition~\ref{prop:dim_and_Y} reduces our main theorem (Theorem~\ref{thm:main})
to the following proposition:
\begin{prop}\label{prop:Yln}
 For $l\in\Z_{\ge1}$ and $n\in\{0,\ldots,l\}$, we have
 \[
  \#Y_{l,n}
  =\frac{1}{l}\sum_{m\mid l}
   \varphi\biggl(\frac{l}{m}\biggr)(2^m-L_m^n).
 \]
\end{prop}

We first invoke the Cauchy-Frobenius lemma:

\begin{prop}[Cauchy-Frobenius lemma]
 If a finite group $G$ acts on a finite set $X$, then we have
 \[
  \#(X/G)
  =\frac{1}{\#G}\sum_{g\in G}\#\{x\in X\mid gx=x\}.
 \]
\end{prop}

\begin{lem}\label{lem:Yln_1}
 For $l\in\Z_{\ge1}$ and $n\in\{0,\ldots,l\}$, we have
 \[
  \#Y_{l,n}
  =\frac{1}{l}\sum_{m\mid l}
   \varphi\biggl(\frac{l}{m}\biggr)
   \#\{\ve{u}\in X_{l,n}\mid m\ve{u}=\ve{u}\}.
 \]
\end{lem}

\begin{proof}
 Applying the Cauchy-Frobenius lemma with $G=\Z/l\Z$ and $X=X_{l,n}$ gives
 \begin{align*}
  \#Y_{l,n}
  &=\frac{1}{l}\sum_{j\in\Z/l\Z}
    \#\{\ve{u}\in X_{l,n}\mid j\ve{u}=\ve{u}\}\\
  &=\frac{1}{l}\sum_{j\in\Z/l\Z}
    \#\{\ve{u}\in X_{l,n}\mid \gcd(j,l)\ve{u}=\ve{u}\}\\
  &=\frac{1}{l}\sum_{m\mid l}
    \bigl(\#\{j\in\Z/l\Z\mid\gcd(j,l)=m\}\cdot
    \#\{\ve{u}\in X_{l,n}\mid m\ve{u}=\ve{u}\}\bigr)\\
  &=\frac{1}{l}\sum_{m\mid l}
    \varphi\biggl(\frac{l}{m}\biggr)
    \#\{\ve{u}\in X_{l,n}\mid m\ve{u}=\ve{u}\}.\qedhere
 \end{align*}
\end{proof}

\begin{lem}\label{lem:Yln_2}
 If $l\in\Z_{\ge1}$, $n\in\{0,\ldots,l\}$, and $m$ is a positive divisor of $l$,
 then
 \[
  \#\{\ve{u}\in X_{l,n}\mid m\ve{u}=\ve{u}\}=2^m-L_m^n.
 \]
\end{lem}

\begin{proof}
 The map $f_{m,l}\colon X_m\to X_l$ defined by
 \[
  f_{m,l}(\ve{v})=(\underbrace{\ve{v},\ldots,\ve{v}}_{l/m})
 \]
 is injective and has image $\{\ve{u}\in X_l\mid m\ve{u}=\ve{u}\}$.
 Therefore it suffices to show that
 \[
  \#f_{m,l}^{-1}(X_{l,n})=2^m-L_m^n.
 \]
 Observe that for each $(m,n)\in\Z_{\ge1}\times\Z_{\ge0}$,
 the set $f_{m,l}^{-1}(X_{l,n})$ is the same for all multiples $l$ of $m$
 with $l\ge n$.
 Therefore, we may put $Z_{m,n}=X_m\setminus f_{m,l}^{-1}(X_{l,n})=f_{m,l}^{-1}(X_{l,n}^c)$,
 aiming to show that $\#Z_{m,n}=L_m^n$ for all $m\in\Z_{\ge1}$ and $n\in\Z_{\ge0}$.

 If $m\le n$, then
 \[
  \#Z_{m,n}
  =\#\{\ve{v}\in X_m\mid\text{at least one component of $\ve{v}$ is $y$}\}
  =2^m-1=L_m^n.
 \]
 If $m=n+1$, then
 \begin{align*}
  \#Z_{m,n}&=\#\{\ve{v}\in X_m\mid\text{at least two components of $\ve{v}$ are $y$}\}\\
  &=2^m-m-1=2^{n+1}-n-2=L_m^n
 \end{align*}
 by Lemma~\ref{lem:n-step} (1).
 Suppose that $m\ge n+2$.
 If we put
 \[
  Z_{m,n,i}=\{(v_1,\ldots,v_m)\in Z_{m,n}\mid v_1=\cdots=v_{i-1}=z,\;v_i=y\}
 \]
 for $i=1,\ldots,n$, then $\{Z_{m,n,i}\mid i=1,\ldots,n\}$
 is a partition of $Z_{m,n}$.
 Set
 \begin{align*}
  Z_{m,n,i}^y&=\{(v_1,\ldots,v_m)\in Z_{m,n,i}\mid v_{i+1}=y\},\\
  Z_{m,n,i}^z&=\{(v_1,\ldots,v_m)\in Z_{m,n,i}\mid v_{i+1}=z\}.
 \end{align*}
 Then removing the $(i+1)$st component gives a bijection from $Z_{m,n,i}^y$ to $Z_{m-1,n,i}$
 and an injection from $Z_{m,n,i}^z$ to $Z_{m-1,n,i}$ with image $Z_{m-1,n,i}\setminus Z_{m-1,n,i}^{z\cdots zy}$,
 where
 \[
  Z_{m-1,n,i}^{z\cdots zy}=\{(v_1,\ldots,v_{m-1})\in Z_{m-1,n,i}\mid v_{i+1}=\cdots=v_{i+n-1}=z,\;v_{i+n}=y\}.
 \]
 Moreover removing components from the $(i+1)$st to the $(i+n)$th gives a bijection from
 $Z_{m-1,n,i}^{z\cdots zy}$ to $Z_{m-n-1,n,i}$.
 It follows that
 \begin{align*}
  \#Z_{m,n,i}&=\#Z_{m,n,i}^y+\#Z_{m,n,i}^z
  =\#Z_{m-1,n,i}+\#(Z_{m-1,n,i}\setminus Z_{m-1,n,i}^{z\cdots zy})\\
  &=2\#Z_{m-1,n,i}-\#Z_{m-n-1,n,i}.
 \end{align*}
 Summing up for $i=1,\ldots,n$ gives $\#Z_{m,n}=2\#Z_{m-1,n}-\#Z_{m-n-1,n}$.
 Hence the lemma follows from Lemma~\ref{lem:n-step} (2).
\end{proof}

Lemmas~\ref{lem:Yln_1} and \ref{lem:Yln_2} imply Proposition~\ref{prop:Yln},
thereby establishing our main theorem.

\section{Multiple zeta-star values}
\subsection{Multiple zeta-star values and Hoffman's algebra}
The \emph{multiple zeta-star values} (MZSVs for short) are defined as
the MZVs with equality allowed in the index of summation:
\begin{defn}
 For $\ve{k}=(k_1,\ldots,k_l)\in I^0$, the \emph{multiple zeta-star value}
 $\zeta^{\star}(k_1,\ldots,k_l)$ is a real number defined by
 \[
  \zeta^{\star}(\ve{k})=\zeta^{\star}(k_1,\ldots,k_l)
  =\sum_{n_1\ge\cdots\ge n_l\ge1}\frac{1}{n_1^{k_1}\cdots n_l^{k_l}}.
 \]
\end{defn}

\begin{ex}
 We have
 \begin{align*}
  \zeta^{\star}(4,2,1)
  &=\sum_{n_1\ge n_2\ge n_3\ge1}\frac{1}{n_1^4n_2^2n_3}\\
  &=\sum_{n_1>n_2>n_3\ge1}\frac{1}{n_1^4n_2^2n_3}
    +\sum_{n_1=n_2>n_3\ge1}\frac{1}{n_1^4n_2^2n_3}
    +\sum_{n_1>n_2=n_3\ge1}\frac{1}{n_1^4n_2^2n_3}
    +\sum_{n_1=n_2=n_3\ge1}\frac{1}{n_1^4n_2^2n_3}\\
  &=\sum_{n_1>n_2>n_3\ge1}\frac{1}{n_1^4n_2^2n_3}
    +\sum_{n_1>n_3\ge1}\frac{1}{n_1^6n_3}
    +\sum_{n_1>n_2\ge1}\frac{1}{n_1^4n_2^3}
    +\sum_{n_1\ge1}\frac{1}{n_1^7}\\
  &=\zeta(4,2,1)+\zeta(6,1)+\zeta(4,3)+\zeta(7).
 \end{align*}
\end{ex}

As the above example indicates,
each MZSV can be expressed as a $\Z$-linear combination of MZVs.
Hoffman's algebra is useful for describing this relationship between MZVs and MZSVs.

\begin{defn}
 Define a $\Q$-linear map $\bar{Z}\colon\H^0\to\R$ by
 setting $\bar{Z}(1)=1$ and $\bar{Z}(z_{k_1}\cdots z_{k_l})=\zeta^{\star}(\ve{k})$
 for $\ve{k}=(k_1,\ldots,k_l)\in I^0$.
\end{defn}

\begin{defn}
 Let $\gamma$ denote the algebra automorphism on $\H$
 satisfying $\gamma(x)=x$ and $\gamma(y)=z$.
 Define a $\Q$-linear transformation $d$ on $\Q+\H y$ by setting
 $d(1)=1$ and $d(wy)=\gamma(w)y$ for $w\in\H$.
\end{defn}

\begin{prop}
 We have $\bar{Z}=Z\circ d\colon\H^0\to\R$.
\end{prop}

\begin{proof}
 Easy and well known.
\end{proof}

\subsection{Cyclic sum formula for multiple zeta-star values}
Ohno and the third author~\cite{OW} proved the following analog of the CSF for MZSVs:
\begin{thm}[Cyclic sum formula for multiple zeta-star values]\label{thm:CSF_MZSVs}
 If $(k_1,\ldots,k_l)\in\check{I}^1$, then
 \[
  \sum_{j=1}^{l}\sum_{i=1}^{k_j-1}\zeta^{\star}(k_j-i+1,k_{j+1},\ldots,k_{j+l-1},i)
  =(k_1+\cdots+k_l)\zeta(k_1+\cdots+k_l+1).
 \]
\end{thm}

The second and third authors~\cite{TW}
defined operators $\bar{\rho}_n$ for MZSVs as well as $\rho_n$ for MZVs:
\begin{defn}
 Let $n\in\Z_{\ge0}$.
 Define a $\Q$-linear map $\bar{\C}_n\colon\H\to\H^{\otimes(n+2)}$
 by setting $\bar{\C}_n(1)=0$, $\bar{\C}_n(x)=x\otimes y^{\otimes(n+1)}$,
 $\bar{\C}_n(y)=-x\otimes y^{\otimes(n+1)}$, and
 \[
  \bar{\C}_n(ww')=\bar{\C}_n(w)\diamond\gamma^{-1}(w')+\gamma^{-1}(w)\diamond\bar{\C}_n(w').
 \]
 Write $\bar{\rho}_n=M_n\circ\bar{\C}_n$.
\end{defn}

\begin{prop}\label{prop:diagram}
 We have the following commutative diagram:
 \[
  \xymatrix{
   &x\H\otimes\H^{\otimes n}\otimes\H y
    \ar[dd]^{\gamma\otimes\gamma^{\otimes n}\otimes d}_{\cong}
    \ar[r]_-{M_n}&\H^0\ar[dd]_{d}^{\cong}\ar[rd]_{\bar{Z}}&\\
   \H\ar[ru]_{\bar{\C}_n}\ar[rd]^{\C_n}\ar@/^4pc/[rru]^{\bar{\rho}_n}\ar@/_4pc/[rrd]_{\rho_n}&&&\R.\\
   &x\H\otimes\H^{\otimes n}\otimes\H y\ar[r]^-{M_n}&\H^0\ar[ru]^{Z}&
  }
 \]
\end{prop}

\begin{proof}
 Straightforward.
\end{proof}

\begin{prop}\label{prop:tedious}
 If $k_1,\ldots,k_l\in\Z_{\ge1}$, then
 \[
  \sum_{j=1}^{l}(\bar{\rho}_0\circ\alpha\circ d)(z_{k_j}\cdots z_{k_{j+l-1}})
  =\sum_{j=1}^{l}\sum_{i=1}^{k_j-1}z_{k_j-i+1}z_{k_{j+1}}\cdots z_{k_{j+l-1}}z_i
   -(k_1+\cdots+k_l)z_{k_1+\cdots+k_l+1}.
 \]
 Here $\alpha\colon\H y\to\H y$ denotes the $\Q$-linear map
 representing the division by depth, i.e.,
 \[
  \alpha(z_{k_1'}\cdots z_{k_{l'}'})=\frac{1}{l'}z_{k_1'}\cdots z_{k_{l'}'}
 \]
 for all monomials $z_{k_1'}\cdots z_{k_{l'}'}\in\H y$.
\end{prop}

\begin{proof}
 For simplicity, let $A$ and $B$ respectively denote
 the left- and right-hand sides of the desired identity.
 By the injectivity of $d$, it suffices to show that $d(A)=d(B)$;
 setting $A'=d(A)$ and $B'=d(B)$, we have
 \begin{align*}
  A'&=\sum_{j=1}^{l}(d\circ\bar{\rho}_0\circ\alpha\circ d)(z_{k_j}\cdots z_{k_{j+l-1}})
     =\sum_{j=1}^{l}(\rho_0\circ\alpha\circ d)(z_{k_j}\cdots z_{k_{j+l-1}}),\\
  B'&=\sum_{j=1}^{l}\sum_{i=1}^{k_j-1}d(z_{k_j-i+1}z_{k_{j+1}}\cdots z_{k_{j+l-1}}z_i)
   -(k_1+\cdots+k_l)z_{k_1+\cdots+k_l+1}.
 \end{align*}

 For $a\in\Z_{\ge1}$, let $\delta_a\colon\H y\to\H y$
 be the $\Q$-linear map that extracts the depth-$a$ part, i.e.,
 \[
  \delta_a(z_{k_1'}\cdots z_{k_{l'}'})=
  \begin{cases}
   z_{k_1'}\cdots z_{k_{l'}'}&\text{if $l'=a$};\\
   0&\text{otherwise}
  \end{cases}
 \]
 for all monomials $z_{k_1'}\cdots z_{k_{l'}'}\in\H y$.
 Then it is enough to prove that $\delta_a(A')=\delta_a(B')$ for $a=1,\ldots,l+1$.
 Note that
 \begin{align*}
  (\delta_a\circ\rho_0)(w)&=
  \begin{cases}
   (\delta_a\circ\rho_0\circ\delta_a)(w)+(\delta_a\circ\rho_0\circ\delta_{a-1})(w)&\text{if $a\ge2$};\\
   (\delta_a\circ\rho_0\circ\delta_a)(w)&\text{if $a=1$},
  \end{cases}\\
  (\delta_a\circ\alpha)(w)&=\frac{1}{a}\delta_a(w)
 \end{align*}
 for all $w\in\H y$.

 For $a=1$, we have
 \begin{align*}
  \delta_1(A')&=\sum_{j=1}^{l}(\delta_1\circ\rho_0\circ\alpha\circ d)(z_{k_j}\cdots z_{k_{j+l-1}})
  =\sum_{j=1}^{l}(\delta_1\circ\rho_0\circ\delta_1\circ\alpha\circ d)(z_{k_j}\cdots z_{k_{j+l-1}})\\
  &=\sum_{j=1}^{l}(\delta_1\circ\rho_0\circ\delta_1\circ d)(z_{k_j}\cdots z_{k_{j+l-1}})
   =\sum_{j=1}^{l}(\delta_1\circ\rho_0)(z_{k_j+\cdots+k_{j+l-1}})\\
  &=\sum_{j=1}^{l}(\delta_1\circ\rho_0)(z_{k_1+\cdots+k_l})
   =\sum_{j=1}^{l}(-z_{k_1+\cdots+k_l+1})=-lz_{k_1+\cdots+k_l+1},\\
  \delta_1(B')
  &=\sum_{j=1}^{l}\sum_{i=1}^{k_j-1}(\delta_1\circ d)(z_{k_j-i+1}z_{k_{j+1}}\cdots z_{k_{j+l-1}}z_i)
   -(k_1+\cdots+k_l)z_{k_1+\cdots+k_l+1}\\
  &=\sum_{j=1}^{l}\sum_{i=1}^{k_j-1}z_{(k_j-i+1)+k_{j+1}+\cdots+k_{j+l-1}+i}
   -(k_1+\cdots+k_l)z_{k_1+\cdots+k_l+1}\\
  &=\sum_{j=1}^{l}(k_j-1)z_{k_1+\cdots+k_l+1}-(k_1+\cdots+k_l)z_{k_1+\cdots+k_l+1}\\
  &=-lz_{k_1+\cdots+k_l+1},
 \end{align*}
 as required.

 For $a=l+1$, we have
 \begin{align*}
  \delta_{l+1}(A')&=\sum_{j=1}^{l}(\delta_{l+1}\circ\rho_0\circ\alpha\circ d)(z_{k_j}\cdots z_{k_{j+l-1}})\\
  &=\sum_{j=1}^{l}(\delta_{l+1}\circ\rho_0\circ\delta_{l+1}\circ\alpha\circ d)(z_{k_j}\cdots z_{k_{j+l-1}})
   +\sum_{j=1}^{l}(\delta_{l+1}\circ\rho_0\circ\delta_l\circ\alpha\circ d)(z_{k_j}\cdots z_{k_{j+l-1}})\\
  &=\frac{1}{l+1}\sum_{j=1}^{l}(\delta_{l+1}\circ\rho_0\circ\delta_{l+1}\circ d)(z_{k_j}\cdots z_{k_{j+l-1}})
   +\frac{1}{l}\sum_{j=1}^{l}(\delta_{l+1}\circ\rho_0\circ\delta_l\circ d)(z_{k_j}\cdots z_{k_{j+l-1}})\\
  &=\frac{1}{l}\sum_{j=1}^{l}(\delta_{l+1}\circ\rho_0)(z_{k_j}\cdots z_{k_{j+l-1}})
  =\frac{1}{l}\sum_{j=1}^{l}\sum_{j'=1}^{l}\sum_{i=1}^{k_{j+j'-1}-1}
   z_{k_{j+j'-1}-i+1}z_{k_{j+j'}}\cdots z_{k_{j+j'+l-2}}z_i\\
  &=\sum_{j=1}^{l}\sum_{i=1}^{k_j-1}z_{k_j-i+1}z_{k_{j+1}}\cdots z_{k_{j+l-1}}z_i
   =\delta_{l+1}(B'),
 \end{align*}
 as required.

 Now let $2\le a\le l$.
 We first compute $\delta_a(A')$.
 It is easy to see that $\delta_a(A')=P+Q$, where
 \begin{align*}
  P&=\frac{1}{a}\sum_{j=1}^{l}
     (\delta_a\circ\rho_0\circ\delta_a\circ d)(z_{k_j}\cdots z_{k_{j+l-1}}),\\
  Q&=\frac{1}{a-1}\sum_{j=1}^{l}
     (\delta_a\circ\rho_0\circ\delta_{a-1}\circ d)(z_{k_j}\cdots z_{k_{j+l-1}}).
 \end{align*}
 For $b\in\Z_{\ge1}$, set
 \[
  M_b=I_{l,b}^1=\{\ve{m}=(m_1,\ldots,m_b)\in\Z_{\ge1}^b\mid m_1+\cdots+m_b=l\}.
 \]
 For each $\ve{m}=(m_1,\ldots,m_b)\in M_b$,
 we extend the indices of $m_p$ to all $p\in\Z$ by declaring $m_p=m_{p'}$
 whenever $p\equiv p'\pmod{b}$.
 For $j=1,\ldots,l$, $\ve{m}=(m_1,\ldots,m_b)\in M_b$, and $p=1,\ldots,b$, write
 \[
  k_{j,\ve{m},p}=\sum_{i=j+m_1+\cdots+m_{p-1}}^{j+m_1+\cdots+m_p-1}k_i
  =k_{j+m_1+\cdots+m_{p-1}}+\cdots+k_{j+m_1+\cdots+m_p-1}.
 \]
 Then we have
 \begin{align*}
  P&=\frac{1}{a}\sum_{j=1}^{l}(\delta_a\circ\rho_0)\Biggl(
     \sum_{\ve{m}\in M_a}z_{k_{j,\ve{m},1}}\cdots z_{k_{j,\ve{m},a}}\Biggr)\\
   &=-\frac{1}{a}\sum_{\substack{j=1,\ldots,l\\\ve{m}\in M_a}}\sum_{p=1}^{a}
      z_{k_{j,\ve{m},p}+1}z_{k_{j,\ve{m},p+1}}\cdots z_{k_{j,\ve{m},p+a-1}}.
 \end{align*}
 Note here that
 \[
  \sum_{\substack{j=1,\ldots,l\\\ve{m}\in M_a}}
  z_{k_{j,\ve{m},p}+1}z_{k_{j,\ve{m},p+1}}\cdots z_{k_{j,\ve{m},p+a-1}}
 \]
 does not depend on $p=1,\ldots,a$,
 because the bijection from $\{1,\ldots,l\}\times M_a$ to itself defined by
 \[
  (j,\ve{m})=\bigl(j,(m_1,\ldots,m_a)\bigr)
  \mapsto(j',\ve{m}')=\bigl(j+m_1+\cdots+m_{p-1},(m_p,\ldots,m_{p+a-1})\bigr)
 \]
 has the property that
 \[
  z_{k_{j,\ve{m},p}+1}z_{k_{j,\ve{m},p+1}}\cdots z_{k_{j,\ve{m},p+a-1}}
  =z_{k_{j',\ve{m}',1}+1}z_{k_{j',\ve{m}',2}}\cdots z_{k_{j',\ve{m}',a}}.
 \]
 It follows that
 \[
  P=-\sum_{\substack{j=1,\ldots,l\\\ve{m}\in M_a}}
     z_{k_{j,\ve{m},a}+1}z_{k_{j,\ve{m},1}}\cdots z_{k_{j,\ve{m},a-1}}.
 \]
 Similar reasoning shows that
 \begin{align*}
  Q&=\frac{1}{a-1}\sum_{j=1}^{l}(\delta_a\circ\rho_0)\Biggl(\sum_{\ve{m}\in M_{a-1}}
     z_{k_{j,\ve{m},1}}\cdots z_{k_{j,\ve{m},a-1}}\Biggr)\\
   &=\frac{1}{a-1}\sum_{\substack{j=1,\ldots,l\\\ve{m}\in M_{a-1}}}
     \sum_{p=1}^{a-1}\sum_{i=1}^{k_{j,\ve{m},p}-1}
     z_{k_{j,\ve{m},p}-i+1}z_{k_{j,\ve{m},p+1}}\cdots z_{k_{j,\ve{m},p+a-2}}z_i\\
   &=\sum_{\substack{j=1,\ldots,l\\\ve{m}\in M_{a-1}}}\sum_{i=1}^{k_{j,\ve{m},a-1}-1}
      z_{k_{j,\ve{m},a-1}-i+1}z_{k_{j,\ve{m},1}}\cdots z_{k_{j,\ve{m},a-2}}z_i.
 \end{align*}

 We next compute $\delta_a(B')$.
 Observe that $\delta_a(B')=R+S$, where
 \begin{align*}
  R&=\sum_{\substack{j=1,\ldots,l\\\ve{m}\in M_{a-1}}}\sum_{i=1}^{k_j-1}
  z_{k_{j,\ve{m},1}-i+1}z_{k_{j,\ve{m},2}}\cdots z_{k_{j,\ve{m},a-1}}z_i,\\
  S&=\sum_{\substack{j=1,\ldots,l\\\ve{m}\in M_a}}\sum_{i=1}^{k_j-1}
  z_{k_{j,\ve{m},1}-i+1}z_{k_{j,\ve{m},2}}\cdots z_{k_{j,\ve{m},a-1}}
  z_{k_{j,\ve{m},a}+i}.
 \end{align*}
 Since the bijection from $\{1,\ldots,l\}\times M_{a-1}$ to itself defined by
 \[
  (j,\ve{m})=\bigl(j,(m_1,\ldots,m_{a-1})\bigr)
  \mapsto(j',\ve{m}')=\bigl(j+m_1,(m_2,\ldots,m_{a-1},m_1)\bigr)
 \]
 has the property that
 \[
  \sum_{i=1}^{k_j-1}
  z_{k_{j,\ve{m},1}-i+1}z_{k_{j,\ve{m},2}}\cdots z_{k_{j,\ve{m},a-1}}z_i
  =\sum_{i=1}^{k_{j'+m_1'+\cdots+m_{a-2}'}-1}
  z_{k_{j',\ve{m}',a-1}-i+1}z_{k_{j',\ve{m}',1}}
  \cdots z_{k_{j',\ve{m}',a-2}}z_i,
 \]
 we have
 \[
  R=\sum_{\substack{j=1,\ldots,l\\\ve{m}=(m_1,\ldots,m_{a-1})\in M_{a-1}}}
  \sum_{i=1}^{k_{j+m_1+\cdots+m_{a-2}}-1}
  z_{k_{j,\ve{m},a-1}-i+1}z_{k_{j,\ve{m},1}}\cdots z_{k_{j,\ve{m},a-2}}z_i.
 \]
 Similar reasoning shows that
 \[
  S=\sum_{\substack{j=1,\ldots,l\\\ve{m}=(m_1,\ldots,m_a)\in M_a}}
  \sum_{i=1}^{k_{j+m_1+\cdots+m_{a-1}}-1}
  z_{k_{j,\ve{m},a}-i+1}z_{k_{j,\ve{m},1}}\cdots z_{k_{j,\ve{m},a-2}}
  z_{k_{j,\ve{m},a-1}+i}.
 \]

 What needs to be shown is that $P+Q=R+S$.
 Note that
 \begin{align*}
  S-P
  &=\sum_{\substack{j=1,\ldots,l\\\ve{m}=(m_1,\ldots,m_a)\in M_a}}
  \sum_{i=0}^{k_{j+m_1+\cdots+m_{a-1}}-1}
  z_{k_{j,\ve{m},a}-i+1}z_{k_{j,\ve{m},1}}\cdots z_{k_{j,\ve{m},a-2}}
  z_{k_{j,\ve{m},a-1}+i}\\
  &=\sum_{\substack{j=1,\ldots,l\\\ve{m}=(m_1,\ldots,m_a)\in M_a}}
  \sum_{i=k_{j,\ve{m},a-1}}^{k_{j,\ve{m},a-1}+k_{j+m_1+\cdots+m_{a-1}}-1}
  z_{k_{j,\ve{m},a-1}+k_{j,\ve{m},a}-i+1}
  z_{k_{j,\ve{m},1}}\cdots z_{k_{j,\ve{m},a-2}}z_i.
 \end{align*}
 Fix $j=1,\ldots,l$ and consider the map $\psi\colon M_a\to M_{a-1}$ defined by
 \[
  \psi(m_1,\ldots,m_a)=(m_1,\ldots,m_{a-2},m_{a-1}+m_a).
 \]
 If $\psi(\ve{m})=\ve{m}'$, then
 $k_{j,\ve{m},a-1}+k_{j,\ve{m},a}=k_{j,\ve{m}',a-1}$ and
 $k_{j,\ve{m},p}=k_{j,\ve{m}',p}$ for $p=1,\ldots,a-2$.
 Moreover, for each $\ve{m}'=(m_1',\ldots,m_{a-1}')\in M_{a-1}$, the sets
 \[
  \{i\in\Z\mid k_{j,\ve{m},a-1}\le i\le k_{j,\ve{m},a-1}+k_{j+m_1+\cdots+m_{a-1}}-1\}
 \]
 for $\ve{m}=(m_1,\ldots,m_a)\in\psi^{-1}(\ve{m}')$
 are disjoint with union
 \[
  \{i\in\Z\mid k_{j+m_1'+\cdots+m_{a-2}'}\le i\le k_{j,\ve{m}',a-1}-1\}.
 \]
 It follows that
 \[
  S-P=\sum_{\substack{j=1,\ldots,l\\\ve{m}'=(m_1',\ldots,m_{a-1}')\in M_{a-1}}}
  \sum_{i=k_{j+m_1'+\cdots+m_{a-2}'}}^{k_{j,\ve{m}',a-1}-1}
  z_{k_{j,\ve{m}',a-1}-i+1}z_{k_{j,\ve{m}',1}}\cdots z_{k_{j,\ve{m}',a-2}}z_i.
 \]
 We therefore conclude that
 \[
  S-P+R=\sum_{\substack{j=1,\ldots,l\\\ve{m}\in M_{a-1}}}
  \sum_{i=1}^{k_{j,\ve{m},a-1}-1}
  z_{k_{j,\ve{m},a-1}-i+1}z_{k_{j,\ve{m},1}}\cdots z_{k_{j,\ve{m},a-2}}z_i
  =Q,
 \]
 as required.
\end{proof}

\begin{prop}\label{prop:rho-bar-ker}
 We have $\bar{\rho}_0(\check{\H}^1)\subset\Ker\bar{Z}$
 and $\bar{\rho}_0(\check{\H}_{k-1}^1)\subset\Ker\bar{Z}\cap\H_k^0$ for $k\in\Z_{\ge2}$.
\end{prop}

\begin{proof}
 If $w\in\check{\H}^1$, then Proposition~\ref{prop:diagram},
 Proposition~\ref{prop:rho0=rho}, and Theorem~\ref{thm:CSF} show that
 \[
  \bar{Z}\bigl(\bar{\rho}_0(w)\bigr)
  =Z\bigl(\rho_0(w)\bigr)
  =Z\bigl(\rho(w)\bigr)
  =0,
 \]
 from which the first assertion follows.
 The second assertion is now obvious.
\end{proof}

\begin{proof}[Proof of Theorem~\ref{thm:CSF_MZSVs}]
 Immediate from Propositions~\ref{prop:tedious} and \ref{prop:rho-bar-ker}.
\end{proof}

\begin{prop}
 We have
 \[
  \{0\}=\bar{\rho}_{k-2}(\check{\H}_1^1)\subset\bar{\rho}_{k-3}(\check{\H}_2^1)
  \subset\cdots\subset\bar{\rho}_0(\check{\H}_{k-1}^1)\subset\Ker\bar{Z}\cap\H_k^0
 \]
 for all $k\in\Z_{\ge2}$.
\end{prop}

\begin{proof}
 Observe that $\bar{\rho}_{n+1}(w)=\bar{\rho}_n(zw)$ for all $n\in\Z_{\ge0}$ and $w\in\H$;
 indeed, we have
 \[
  \bar{\rho}_{n+1}(w)=d^{-1}\bigl(\rho_{n+1}(w)\bigr)=d^{-1}\bigl(\rho_n(zw)\bigr)=\bar{\rho}_n(zw).
 \]
 This finishes the proof.
\end{proof}

\begin{thm}
 If $n\in\Z_{\ge0}$ and $k\in\Z_{\ge1}$, then
 \[
  \dim\bar{\rho}_n(\check{\H}_k^1)
  =\dim\rho_n(\check{\H}_k^1)
  =\frac{1}{n+k}\sum_{m\mid n+k}\varphi\biggl(\frac{n+k}{m}\biggr)(2^m-L_m^n)-2.
 \]
\end{thm}

\begin{proof}
 The $\Q$-vector spaces $\bar{\rho}_n(\check{\H}_k^1)$ and $\rho_n(\check{\H}_k^1)$
 are isomorphic as given by $d$,
 and the theorem follows from Theorem~\ref{thm:main}.
\end{proof}

\section{Acknowledgements}
The authors are grateful to Professor Yasuo Ohno for valuable comments that inspired
this work, and to Mr.~Hiroki Kondo for carefully reading the manuscript.

\bigskip
\hrule
\bigskip

\noindent 2010 \textit{Mathematics Subject Classification}:
Primary 11M32; Secondary 05A05.

\noindent \textit{Keywords}:
multiple zeta values, multiple zeta-star values, cyclic sum formula, Lucas numbers.

\bigskip
\hrule
\bigskip

\noindent (Concerned with sequences
\seqnum{A000931}, \seqnum{A000079}, \seqnum{A038360}, \seqnum{A052823}, \seqnum{A000012},
\seqnum{A000032}, \seqnum{A000204}, \seqnum{A001644}, \seqnum{A073817}, \seqnum{A001648},
\seqnum{A074048}, \seqnum{A023424}.)
\end{document}